\documentclass[11pt]{article}
\usepackage{amsmath,amsfonts,amssymb,theorem,color,amsrefs}
\newtheorem{theorem}{Theorem}

\newtheorem{lemma}[theorem]{Lemma}

\newcommand{\qed}{\hfill\fbox{}\par\vspace{0.3mm}}
\newenvironment{proof}{{\bf Proof.}} {\hfill\qed}

\numberwithin{equation}{section}
\numberwithin{theorem}{section}
\numberwithin{definition}{section}

\newcommand{\de}{\delta}
\newcommand{\e}{\epsilon}
\newcommand{\ga}{{\gamma}}
\newcommand{\Om}{{\Omega}}
\newcommand{\om}{{\omega}}
\newcommand{\si}{\sigma}
\renewcommand{\th}{\theta}
\newcommand{\la}{\lambda}
\newcommand{\al}{\alpha}

\newcommand{\bke}[1]{\left( #1 \right)}
\newcommand{\bkt}[1]{\left[ #1 \right]}
\newcommand{\bket}[1]{\left\{ #1 \right\}}
\newcommand{\norm}[1]{\| #1 \|}

\newcommand{\R}{{\mathbb R }}
\newcommand{\N}{{\mathbb N}}
\newcommand{\pd}{{\partial}}
\newcommand{\nb}{{\nabla}}

\newcommand{\ph}{{\varphi}}
\newcommand{\lv}{{\bar v}}
\newcommand{\lp}{{\bar p}}
\newcommand{\lom}{{\bar \omega}}
\newcommand{\td}{\tilde}

\renewcommand{\div}{\mathop{\mathrm{div}}}
\newcommand{\curl}{\mathop{\mathrm{curl}}}

\newcommand{\osc}{\mathop{\mathrm{osc}}}

\newcommand{\donothing}[1]{{}}

\begin{document}
\title{Lower bound on the blow-up rate of the axisymmetric \\ Navier-Stokes equations}

\author{Chiun-Chuan Chen\thanks{
National Taiwan University and National Center for Theoretical Sciences, 
Taipei Office, email: chchchen@math.ntu.edu.tw}\;
, \quad Robert M. Strain\thanks{Harvard University, 
email: strain@math.harvard.edu}\;,\quad 
Tai-Peng Tsai\thanks{University of British Columbia, email: ttsai@math.ubc.ca}\;,
\quad Horng-Tzer Yau\thanks{Harvard University, email: htyau@math.harvard.edu}
}

\date{}

\maketitle


\begin{abstract} Consider axisymmetric strong solutions of the incompressible Navier-Stokes equations  in $\R^3$ with non-trivial swirl.  Such solutions are not known to be globally defined, but it is shown in \cite{MR673830} that they could only blow up on the axis of symmetry.
  Let $z$ denote the axis of symmetry and $r$ measure the distance 
to the $z$-axis.   Suppose the solution satisfies the pointwise scale invariant bound
$|v (x,t)| \le C_*{(r^2 -t)^{-1/2}} $ for $-T_0\le t < 0$ and $0<C_*<\infty$ allowed to be large,  we then prove that $v$ is regular at time zero. 
\end{abstract}


\section{Introduction}

The incompressible Navier-Stokes equations in {\it cartesian coordinates} are given by
\begin{equation} \tag{N-S}
\label{nse}
\pd_t v + (v \cdot \nb)v + \nb p = \Delta v, \quad \div v =0.
\end{equation}
The velocity field is $v(x,t)=(v_1, v_2, v_3):\mathbb{R}^3\times [-T_0,0) \to \R^3$  and
$p(x,t):\mathbb{R}^3\times [-T_0,0) \to \R$ is the pressure.  It is a long standing open question
to determine if solutions with large smooth initial data of finite energy
remain regular for all time.

In this paper we consider the special class of solutions which are {\it axisymmetric}.  This means, in {\it  cylindrical coordinates} $r,\th,z$ with
$(x_1,x_2,x_3)=(r\cos \th,r\sin \th,z)$, that the solution is of the form
\begin{equation} \label{v-cyl}
v(x,t) = v_r(r,z,t)e_r + v_\th(r,z,t)e_\th + v_z(r,z,t)e_z.
\end{equation}
The components $v_r,v_\th,v_z$ do not depend upon $\th$ and the basis
 vectors $e_r,e_\th,e_z$ are 
\[
e_r = \left(\frac {x_1}r,\frac {x_2}r,0\right),\quad
e_\th = \left(-\frac {x_2}r,\frac {x_1}r,0\right),\quad
e_z = (0,0,1).
\]
The main result of our paper 
 shows that axisymmetric solutions must blow up faster
than the scale invariant rate which appears in \eqref{assumption} below.

For $R>0$ define $B(x_0,R)\subset \R^{3}$ as the ball of radius $R$
centered at $x_0$.  The parabolic cylinder is  $Q(X_0,R) = B(x_0,R)\times (t_0-R^2,t_0)\subset
\R^{3+1}$  centered at $X_0=(x_0,t_0)$. If the center is
the origin we use the abbreviations  $B_R=B(0,R)$ and $Q_R=Q(0,R)$.

\begin{theorem}\label{mainthm}
Let $(v,p)$ be an axisymmetric solution of the Navier-Stokes equations
\eqref{nse}  in $D=\R^3 \times (-T_0,0)$ for which
$v(x,t)$ is smooth in $x$ and H\"older continuous in $t$. Suppose the pressure satisfies $p \in
L^{5/3}(D)$ and $v$ is pointwise bounded as
\begin{equation}\label{assumption}
|v (x,t)| \le C_*{(r^2 -t)^{-1/2}}, \quad (x,t) \in D.
\end{equation}
The constant $C_*<\infty$ is allowed to be large. Then $v
\in L^\infty(B_R \times [-T_0,0])$ for any $R>0$.
\end{theorem}

We remark that the exponent $5/3$ for the norm of $p$ can be replaced.  However,
it is the natural exponent occurring in the existence theory for weak solutions,  see e.g. \cite{MR0454426}, 
\cite{MR673830}.

Recall the natural scaling of Navier-Stokes equations:  If $(v,p)$ is a solution to \eqref{nse}, then for any $\la>0$
the following rescaled pair is also a solution:
\begin{equation}
\label{rescale}
v^\la(x, t) = \la v(\la x, \la^2 t),\quad
p^\la(x,t) = \la^2 p (\la x, \la^2 t).
\end{equation}
Suppose a
solution $v(x,t)$ of the Navier-Stokes equations  blows up at  $X_0=(x_0,t_0)$.  Leray
\cite{JFM60.0726.05} proved that the blow up rate in time is at least
\[
\norm{v(\cdot,t)}_{L^\infty_x} \ge C{(t_0-t)^{-1/2}}.
\]
Caffarelli, Kohn, and Nirenberg \cite{MR673830}  showed that for such a
blow-up solution the average of $|v|$ over $Q(X_0,R)$ satisfies 
\[
\bigg( \frac 1{|Q_R|} \int_{Q(X_0,R)} |v|^3 + |p|^{3/2} dx dt \bigg)^{1/3} \ge \frac CR.
\]
See also \cite{MR1488514,MR1397564,Vasseur}.
Thus, the natural rate for blow-up is at least
\begin{equation}\label{eq1-3}
|v(x,t)| \sim \frac {O(1)}{  [(x_0-x)^2 + t_0 -t]^{1/2} }.
\end{equation}
Both this and the rate \eqref{assumption} are invariant under the natural scaling  \eqref{rescale}.

The Serrin type criteria \cite{MR0136885,MR0236541,MR0316915, MR833416,MR716283,MR933230,GKT}
states that $v$ is regular if it satisfies
\begin{equation}\label{serrin}
\norm{v }_{L^s_tL^q_x(Q_1)}< \infty,
\quad \frac{3}{q}+\frac{2}{s}\le 1, \ s,q \in (2,\infty), \quad\text{or} \quad(s,q)=(2,\infty).
\end{equation}
Above,  for a domain $\Om \subset \R^3$, we use the definition
\[
\norm{v }_{L^s_tL^q_x(\Om \times (t_1,t_2))} := \norm{
\norm{v(x,t)}_{L^q_x(\Om)}}_{L^s_t(t_1,t_2)}.
\]
For any $X_0=(x_0,t_0)\in Q_1$, \eqref{serrin} implies the following local smallness of $v$:
\begin{equation}\label{ereg}
\lim_{R\downarrow 0}\norm{v }_{L^s_tL^q_x(Q(X_0,R))} =0.
\end{equation}
Therefore \eqref{serrin} is a
 so-called $\e$-regularity criterion since it implies that the norm is locally small.    For $(q, s)=(3, \infty)$, \eqref{ereg}
does not follow from \eqref{serrin}.  Hence the $(q, s)=(3, \infty)$ end point regularity
criterion \eqref{serrin} proved in \cite{MR1972005,MR2005639} is not an $\e$-regularity type theory.

However these criteria do not rule out blow-up with the natural
scaling rate \eqref{eq1-3}. 
It is a fundamental problem in the study
of the incompressible Navier-Stokes equations to determine if  solutions to \eqref{nse} with the following scale invariant bound are regular
\begin{equation}\label{nb}
|v(x,t)| \le \frac C{  [(x_0-x)^2 + t_0 -t]^{1/2} }.
\end{equation}
If a self-similar solution satisfies  this bound
then it is known to be zero \cite{MR1643650} (the self-similar solution from \cite{MR1397564} belongs to $L_t^\infty L_x^3$).

Theorem \ref{mainthm} rules out singular axisymmetric solutions satisfying the bound \eqref{nb}. In fact
\eqref{assumption} is considerably weaker than \eqref{nb} and 
is also not a borderline
case of the Serrin type criterion. For example  \eqref{assumption} implies that $v \in L^q(Q_1)$ for $q<4$, but not for
$q\in [4,5)$. The borderline of the Serrin type criterion, on the other hand,  is $v \in L^5(Q_1)$.

We now recall the previous results on the regularity of axisymmetric
solutions to the Navier-Stokes equations.  Global in time
regularity was first proved under the {\it no
swirl} assumption,  $v_\th=0$, independently by
Ukhovskii-Yudovich \cite{MR0239293} and Ladyzhenskaya \cite{MR0241833}. See
\cite{MR1718156} for a refined proof and \cite{MR2083793} for similar results in
the half space setting. 

When the swirl component $v_\th$ is not
assumed to be trivial, global regularity is unknown.  But it follows from
the partial regularity theory of \cite{MR673830} that singular points can
only lie on the axis of symmetry.  Any off axis symmetry would imply a whole circle of singular points, which contradicts \cite{MR673830}.
Neustupa-Pokorn\'y \cite{MR1814224,MR1844284} proved
regularity assuming the zero dimensional condition $v_r \in L^s_tL^q_x$ with $3/q+2/s
=1$, $3<q\le \infty$. 
Regularity criteria can also be put on the vorticity field $\om =
\curl v$:
\begin{equation}\label{om-cyl}
\om(x,t) = \om_r e_r + \om_\th e_\th + \om_z e_z,   
\end{equation}
where
\[
\om_r = - \pd_z v_\th, \quad \om_\th = \pd_z v_r - \pd_r v_z, 
\quad \om_z = (\pd_r + r^{-1})v_\th.
\]
Chae-Lee \cite{MR1902055} proved regularity assuming finiteness of another
zero-dimensional integral: $\om_\th \in L^s_tL^q_x$ with $3/q+2/s =2$. 
Jiu-Xin  \cite{MR2055842} proved regularity if the sum of the
zero-dimensional {\it scaled} norms $ \int_{Q_R} (R^{-1}|\om_\th|^2 +R^{-3}
|v_\th|^2) dz$ is sufficiently small for $R>0$ small enough.
Recently, Hou-Li \cite{HouLi} constructed a family of global solutions with large initial data.

The main idea of our proof is as follows. The bound \eqref{assumption}
ensures that the first blow up time is no earlier than $t=0$.  For
$t\in (-T_0, 0)$ we show that the swirl component $v_\th$ gains a modicum
of regularity: For some small $\alpha=\alpha(C_*)>0$,  \eqref{assumption} enables us to conclude that
\begin{equation}
|v_\theta(t,r,z)| \le C r^{\alpha-1}.
\label{modicum}
\end{equation}
We prove \eqref{modicum} in Section \ref{sec:holder}.  This estimate breaks the scaling, thereby transforming the problem from order one to $\epsilon$-regularity, which is shown to be sufficient in Section \ref{pfOFmainTHM}.

\section{Proof of main theorem}\label{pfOFmainTHM}

In this section we prove Theorem \ref{mainthm}.  First we show that our solutions are in fact suitable weak solutions.  Then we make use of
 \eqref{modicum}, to establish our main theorem.

\subsection{Suitable weak solution}

We recall from \cite{MR0454426,MR673830,MR1488514} that a {\it suitable weak solution} of
the Navier-Stokes equations in a domain $Q \subset \R^3 \times \R$ is
defined to be a pair $(v,p)$ satisfying
\begin{equation}\label{SWSbound}
v \in L^\infty_t L^2_x(Q), \quad 
\nb v \in L^2(Q), \quad p \in L^{3/2}(Q).
\end{equation}
Further $(v,p)$ solve \eqref{nse} in the sense of distributions and
satisfy the {\it local energy inequality}:
\begin{equation}
 2\int_{Q} |\nabla v|^2 \varphi
\le
\int_{Q} \left\{|v|^2(\partial_t\varphi+\Delta \varphi)
+
 ( |v|^2+2p)v\cdot\nabla\varphi\right\}, \quad \forall \ph \in
C^\infty_c(Q), \ \ph \ge 0.
\label{localE}
\end{equation}
To prove interior regularity, we do not need to specify
the initial or boundary data.

We define a solution $v(x,t)$ to be {\it regular} at a point $X_0$ if $v \in
L^\infty(Q(X_0,R))$ for some $R>0$.  Otherwise $v(x,t)$ is {\it
singular} at $X_0$. We will use the following regularity criterion.

\begin{lemma}\label{th2-1}
Suppose that $(v,p)$ is a suitable weak solution
of \eqref{nse} in $Q(X_0,1)$.  Then there exists an $\e_1>0$ so that $X_0$ is a regular point if
\begin{equation} \label{eq2-4-1}
\limsup_{R \downarrow 0} \frac{1}{R^2} \int_{Q(X_0,R)} |v|^3 \le \e_1.
\end{equation}
\end{lemma}

This regularity criterion, which is a variant of the criterion in \cite{MR673830}, was proven in \cite{MR1685610};  see \cite{GKT}
for more general results. The condition \eqref{eq2-4-1} does
not explicitly involve the pressure, but one does require $p \in L^{3/2}(Q(X_0,1))$ because the pair $(v,p)$ is assumed to be a suitable weak solution.

\subsection{Preliminary estimates}

In this subsection we show that the solution $(v,p)$ in Theorem
\ref{mainthm} is sufficiently integrable to be a suitable weak solution, and we derive estimates depending only upon $C_*$ of \eqref{assumption}.

We  estimate the pressure with weighted singular integral estimates.  We therefore first estimate $v$ in weighted spaces.
Fix $\beta \in (1,5/3)$. For $t \in (-T_0,0)$ by \eqref{assumption} we have
\[
\int_{\R^3} \frac{|v(x,t)|^4 }{|x|^{\beta}} \,dx \le \int_{\R^3} \frac
{1}{|x|^\beta} \frac{C_* rdr dz}{(r^2-t)^2} 
=\int_{|z|\ge 1}+\int_{|z|< 1, r>1}+\int_{|z|< 1, r<1}
= I_1+I_2+I_3.
\]
Each of these integrals can be estimated as follows
\begin{align*}
|I_1|&\le 
\int_{|z|>1} \frac{dz}{|z|^{\beta}} \int_{0}^\infty  \frac{C_* r dr}{(r^2-t)^2}
\le c|t|^{-1},
\\
|I_2|&\le 
  \int_{1}^\infty r^{-\beta} \frac{C_*}{(r^2-t)^2}rdr
\le c,
\\
|I_3|&\le 
\int_{0}^1 (1 + r^{1-\beta})  \frac c{(r^2-t)^2}rdr
\le c|t|^{-(1+\beta)/2} .
\end{align*}
Summing the estimates and using $\beta>1$ we get
\[
\int_{\R^3} \frac 1{|x|^{\beta}} |v(x,t)|^4 \,dx 
\le c + c|t|^{-(1+\beta)/2}.
\]
Define $R_i$'s to be the Riesz transforms: $R_i = \frac{\partial_i}{\sqrt{-\Delta}}$.   We consider the singular integral
\[
\td p(x,t) = \int \sum_{i,j}\pd_i  \pd_j (v_i v_j)(y) \frac
1{4\pi|x-y|} \, dy = \sum_{i,j}R_i R_j(v_i v_j ).
\]
To show that this singular integral is well defined for every $t$, we use the $L^q(\mathbb{R}^3)$-estimates for singular integrals with $A_q$ weight \cite{MR1232192}.  Specifically, we use $q=2$ and the $A_2$ weight function $|x|^{-\beta}$.  We have the estimate
\begin{equation}\label{eq4-1}
\int \frac 1{|x|^{\beta}} |\td p(x,t)|^2 \,dx \le
c\int \frac 1{|x|^{\beta}} |v(x,t)|^4 \,dx \le c+c|t|^{-(1+\beta)/2} .
\end{equation}
Choose $\ga \in (1/2 + 5\beta/6, 3)$.  H\"older's inequality gives us the bound
\[
\int_{|x|>1} \frac  {|\td p(x,t)|^{5/3}}{|x|^{\ga} }dx \le \bke{\int_{|x|>1}
\frac {|\td p(x,t)|^2}{|x|^\beta}dx }^{5/6} \bke{\int_{|x|>1}
|x|^{-(\ga-\frac 56\beta)6} dx}^{ 1/6} < \infty.
\]
We will use these bounds to show that the pressure $p$ can be identified with $\td p$.

Let $h(x,t) = p(x,t) - \td p(x,t)$. Then $h$ is harmonic in $x$, $\Delta_x h(x,t) = 0$, and by assumption
$p(\cdot,t) \in L^{5/3}(\R^3)$ for almost every $t$. For each such $t$
we have
\[
\int_{|x|>1} \frac {|h(x,t)|^{5/3}}{|x|^\ga} dx 
\le c\int_{|x|>1} |p(x,t)|^{5/3}dx 
+ c\int_{|x|>1} \frac {|\td p(x,t)|^{5/3}}{|x|^\ga} dx 
<\infty.
\]
We may thus conclude from using a Liouville theorem that $h(x,t) =0$ for all $x$ if $\ga<3$. 

To see the last
assertion, fix a radial smooth function $\phi(x) \ge 0$ supported in
$2< |x| <4$ satisfying $\int \phi =1$. For any $x\in \R^3$
with $R>|x|$ we have
\[
h(x,t) = \int h(y,t)R^{-3}\phi(x+y/R)\,dy.
\]
This is the mean value theorem for harmonic functions.
Define $A= B_{5R}-B_R$, then 
\[
|h(x,t)| \le cR^{-3} \int_A |h(y,t)| dy\le cR^{-3+(6+3\ga)/5}
\bke{\int_A |y|^{-\ga} |h(y,t)|^{5/3} dy}^{3/5}.
\]
This clearly vanishes as $R \to \infty$.  Thus  $p(x,t)=\td p(x,t)$
for all $x$ and for almost every $t$.

Next we show that $(v,p)$ form a suitable weak solution.
From H{\"o}lder's  inequality,  \eqref{eq4-1} and $\beta<5/3$ we conclude that
\begin{equation}
\int_{Q_{1}} |p(x,t)|^{3/2}\,dxdt \le c \int_{-1}^0   \bke{ 
\int_{B_{1}}  \frac 1{|x|^{\beta}} |p(x,t)|^2 \,dx}^{3/4} dt \le c .
\label{p32}
\end{equation}
The pointwise estimate \eqref{assumption} on $v$ implies
\begin{equation}\label{vbd}
v \in L^s_tL^q_x(Q_1), \quad \frac{1}{q}+\frac{1}{s}> \frac{1}{2}.
\end{equation}
We will use $(s,q)=(3,3)$.  We also see from \eqref{assumption} that $v \in L^4( B_1 \times (-T_0,-\e))$ for any small $\e >0$. Thus the
vector product of \eqref{nse} with $u \ph$ for any $\ph \in
C^\infty_c(Q_{1})$ is integrable in $Q_{1}$ and we can integrate by
parts to get the local energy inequality \eqref{localE} with
$Q=Q_{1}$. In fact we have equality.

Now, for any $R \in (0,1)$ and $t_0 \in (-R^2,0)$, we can choose
a sequence of $\ph$ which converges a.e.~in $Q_R$ to $H(t_0-t)$,
the Heviside function that equals $1$ for $t<t_0$ and $0$ for
$t>t_0$. Since the limit of $\pd_t\ph$ is the negative delta function
in $t$, this gives us the estimate
\begin{equation}\label{SWS-v}
\mathop{\mathrm{ess\, sup}}_{-R^2<t<0}\int_{B_R}|v(x,t)|^2dx
+  \int_{Q_R} |\nb v|^2 \le C _{R} \int _{Q_1} (|v|^3 +
|p|^{3/2}).
\end{equation}
These estimates show that $(v,p)$ is a suitable weak solution of \eqref{nse}
in $Q_{R}$. Note that these bounds depend on $C_*$ of
\eqref{assumption} only, not on 
$\norm{p}_{L^{5/3}(\R^3 \times (-T_0,0))}$.

\subsection{Scaling limit}\label{sec:slimit}

To show Theorem \ref{mainthm}, it suffices to show that every point on
the $z$-axis is regular.  Suppose now a point $x_*=(0,0,x_3)$ on the
$z$-axis is a singular point of $v$. We will derive a
contradiction. Define $X_*=(x_*,0)$. Let $(v^\la,p^\la)$ be rescaled solutions of
\eqref{nse} defined by
\begin{equation}
\label{rescale-1}
v^\la(x,t) = \la v(\la (x-x_*), \la^2 t),\quad
p^\la(x,t) = \la^2 p (\la (x-x_*), \la^2 t).
\end{equation}
By Lemma \ref{th2-1}, there is a sequence $\la_k$, $k\in \N$, so that $\la_k
\to 0$ as $k \to \infty$ and
\begin{equation}
\label{blowup0}
\int_{Q_1} |v^{\la_k}|^3 = \frac{1}{\la_k^2} \int_{Q(X_*,\la_k)} |v|^3
>\e_1.
\end{equation}
We will derive a contradiction to this statement.

For $(v^\la,p^\la)$ with $0<\la<1$, the
pointwise estimate \eqref{assumption} is preserved:
\[
|v^\la(x,t)| \le C_* (r^2-t)^{-1/2}, \quad  (x,t) \in \R^3 \times (-T_0,0).
\]
We also have by rescaling
\[
p^\la(x,t) = \int \sum_{i,j}\pd_i  \pd_j (v_i^\la v_j^\la)(y) \frac
1{4\pi|x-y|} \, dy ,
\]
The argument in the previous subsection provides the uniform
bounds for $q\in (1,4)$: 
\begin{equation}\label{u-bd}
  \int_{Q_1}   |v^\la |^q + |p^\la|^{3/2}  \le C,
\quad
\mathop{\mathrm{ess\, sup}}_{-R^2<t<0}
\int_{B_R}|v^\la(x,t)|^2dx +  \int_{Q_R} |\nb v^\la|^2 \le C.
\end{equation}
Above the bound for $p_\lambda$ follows from \eqref{p32}, the bound for 
$|v^\lambda|^q$ follows from \eqref{assumption}, and the energy bound then follows from \eqref{SWS-v}.

Thus from the sequence $\la_k$ we can extract a subsequence, still
denoted by $\la_k$, so that $(v^{\la_k},p^{\la_k})$ weakly converges
to some limit function $(\lv,\lp)$
\[
v^{\la_k} \rightharpoonup \lv \quad \text{in }L^q(Q_R),\quad
\nb v^{\la_k} \rightharpoonup \nb \lv \quad \text{in }L^2(Q_R),\quad
p^{\la_k} \rightharpoonup \lp \quad \text{in }L^{3/2}(Q_R).
\]
Moreover since $(v^\la,p^\la)$ solves \eqref{nse} with bound
\eqref{u-bd}, we also have the uniform bound
\[
\norm{\pd_t v^\la}_{L^{3/2}((-R^2,0); H^{-2}(B_R))} < C.
\]
We can then apply Theorem 2.1 of \cite[chap.~III]{MR0609732} to conclude that the
$v^{\la_k}$ remain in a compact set of $L^{3/2}(Q_R)$. Therefore (a
further subsequence of) $v^{\la_k}\to \lv$ strongly in
$L^{3/2}(Q_R)$. Since the $v^{\la_k}$ remain bounded in $L^q(Q_R)$
for all $q<4$, we deduce that $v^{\la_k}\to \lv$ strongly in
$L^{q}(Q_R)$ for all $1\le q <4$.

\subsection{The limit solution}

The convergence established at the end of Section \ref{sec:slimit} is sufficient  to conclude that the limit
function $(\lv,\lp)$ is a suitable weak solution of the Navier-Stokes
equations in $Q_R$, as in \cite{MR673830,MR1488514}.  Since $v$ satisfies \eqref{assumption}  so
does $\lv$. Hence $\lv$ is regular at any interior point of $Q_R$,
and $t=0$ is the first time when $\lv(x,t)$ could develop a
singularity.

To gather more information we use axisymmetry.  We will argue in this subsection and the next that the 
 estimate \eqref{modicum} (proven in the Section \ref{sec:holder}) is enough to conclude that our solution is regular.  In particular \eqref{modicum}
tells us that
\[
\int_{Q_R} \left|v^\la_\theta\right| \le C\la^\alpha \to 0
\quad \text{as} \quad \la \downarrow 0.
\] 
Thus the limit $\lv$ has no-swirl, $\lv_\theta =0$.

Let $\lom = \nb \times \lv$ be the vorticity of $\lv$.  The
$\theta$ component of $\lom$, $\lom_\theta =\pd_z \lv_r- \pd_r
\lv_z$, solves
\begin{equation*}
 \left(\pd_t + \bar{b}\cdot \nabla - \Delta + \frac 1{r^2}\right) \bar{\om}_\th 
 - \frac{\bar{v}_r}r \bar{\om}_\th = 0.
\end{equation*}
We have used $\lv_\th=0$.
Above
\begin{equation*} 
\bar{b}= \lv = \lv_r e_r + \lv_z e_z,
\quad
\bar{b}\cdot \nabla =\lv_r \partial_r+\lv_z \partial_z, \quad \div \bar{b} = 0.
\end{equation*}
We record the Laplacian for axisymmetric functions
$$
\Delta
=
 \frac{\partial^2}{\partial r^2}+\frac{1}{r}\frac{\partial}{\partial r}
+\frac{\partial^2}{\partial z^2}.
$$
Next define $\Om = \lom_\theta/r$.  Then $\Om$ solves
\begin{equation}\label{Om2}
\left(\pd_t + \bar{b} \cdot \nb - \Delta - {\frac 2r} \pd_r \right) \Om=0.
\end{equation}
We now derive $L^q$ estimates on $\Om$ using estimates for the Stokes system.

Since $\bar{v}$ satisfies \eqref{assumption}, it also satisfies \eqref{vbd}.
We will use both $(s,q)=(5/2,5)$ and $(s,q)=(5/4,5/4)$ below. 
We rewrite \eqref{nse} as a Stokes system with force
\[
 (\pd_t - \Delta )\lv_i + \nb_i \lp = \pd_j f_{ij}, 
 \quad \div \lv = 0, 
\quad f_{ij} = -\lv_i \lv_j.
\]
By the interior estimates of Stokes system (shown in the Appendix) we have
\[
\norm{ \nb \lv}_{L^{5/4}_tL^{5/2}_x (Q_{5/8})} \le 
C \norm{  \lv}_{ L^{5/2}_tL^{5}_x (Q_{3/4})}^2
+ C \norm{  \lv}_{ L^{5/4}(Q_{3/4})} \le C.
\]
Hence $\Om$ has the bound
\begin{equation}\label{Om3}
\norm{\Om}_{L^{20/19}(Q_{5/8})} \le \norm{ \nb \lv}_{L^{5/4}_tL^{5/2}_x
(Q_{5/8})} \norm{1/r}_{L^{\infty}_tL^{20/11}_x (Q_{5/8})} \le C.
\end{equation}
In Section \ref{sec:wH}, we 
obtain $\Om\in L^\infty$ from  \eqref{Om2}, \eqref{Om3} and a local maximum estimate.  Then in 
Section \ref{sec:regO} we show that this is sufficient 
to conclude Theorem \ref{mainthm}.

\subsection{Energy Estimates}

We derive parabolic De Giorgi type energy estimates for  \eqref{Om2}.  To do this we assume that
\begin{equation*} 
|\bar{b}(r,z,t)|\le C_*/r.
\end{equation*}
This assumption on $\bar{b}$ is substantially weaker than the one from Theorem \ref{mainthm}. 

Consider a  test function $0\le \zeta(x,t)\le 1$ defined on $Q_R$ for which $\zeta=0$ on 
$\pd B_R \times [- R^2,0]$ and $\zeta =1$ on $Q_{\sigma R}$ for $0<\sigma<1$.  
Define $(u)_\pm = \max\{ \pm u, 0\}$ for a scalar function $u$.
Multiply \eqref{Om2} by $p (\Omega-k)_\pm^{p-1} \zeta^2$ for  $1<p \le 2$ and $k\ge 0$ to obtain
\begin{gather*}
\left.\int_{B_R} \zeta^2 (\Omega-k)_\pm^p\right|^0_{- R^2}
+
\frac {4(p-1)}p \int_{Q_R} |\nabla
((\Omega-k)_\pm^{p/2}\zeta)|^2 
\\
=
2\int_{Q_R}  (\Omega-k)_\pm^p 
\left(
\zeta \frac{\partial \zeta}{\partial t}+|\nabla \zeta|^2 + \frac{2-p}p \zeta\Delta \zeta
-2\zeta \frac{\partial_r \zeta}{r} 
+\bar{b}  \cdot \zeta \nabla \zeta\right)
\\
-
 2\left.\int dt\int dz ~ \zeta^2 (\Omega-k)_\pm^p\right|_{r=0}.
\end{gather*}
Notice that the last term has a good sign.

To estimate the term involving $b$ we use  Young's inequality
$$
\int_{\mathbb{R}^3} v_\pm^2 b  \zeta\cdot \nabla \zeta
\le \delta \frac{R^{-1+\e}}{1+\epsilon} \int_{\mathbb{R}^3} v_\pm^2 \zeta^2 |b|^{1+\e} + 
C_\delta\frac{\epsilon R^{-2+ (1+\e)/\e}}{1+\epsilon}
 \int_{\mathbb{R}^3} v_\pm^2  \zeta^2  \left [  \frac {|\nabla \zeta|} \zeta \right ]^{ (1 + \e)/\e}.
$$
This holds  for small $\delta>0$ and $\epsilon>0$ to be chosen.
Further choose $\zeta$ to decay like  $( 1-|x|/R)^n $  near the boundary of $B_R$.
If $n$ is large enough
(depending on $\e$)  we have
$$
C_\delta\frac{\epsilon R^{-2+ (1+\e)/\e}}{1+\epsilon}
 \int_{\mathbb{R}^3} v_\pm^2  \zeta^2  \left [  \frac {|\nabla \zeta|} \zeta \right ]^{ (1 + \e)/\e}
 \le
 C
R^{-2} \int_{B_R}  v_\pm^2.
$$
We also use the H{\"o}lder and Sobolev inequalities to obtain 
\begin{equation*}
\begin{split}
\delta \frac{R^{-1+\e}}{1+\epsilon} \int_{\mathbb{R}^3} v_\pm^2 \zeta^2 |b|^{1+\e}
&\le
\delta\left(R^{(-1+\e)3/2}  \int_{B_R}   |b|^{(1+\e)3/2}\right)^{2/3} \int_{\mathbb{R}^3} |\nabla (v_\pm\zeta ) |^2
\\
&\le 
\delta C\int_{\mathbb{R}^3} |\nabla (v_\pm\zeta ) |^2.
\end{split}
\end{equation*}
The last inequality is satisfied for example if $|b|\le C_*/r$ and
$\epsilon <1/3$.  We conclude
\begin{equation}
\label{b-term}
\begin{split}
\int_{\mathbb{R}^3} v_\pm^2 b  \zeta\cdot \nabla \zeta
&\le 
\delta C\int_{\mathbb{R}^3} |\nabla (v_\pm\zeta ) |^2
+
 C
R^{-2} \int_{B_R}  v_\pm^2.
\end{split}
\end{equation}
The key point which we used here to control the more
singular drift term was to split $b$ from the main part of the term
$v_\pm\zeta$, using the Young and Sobolev inequalities instead of
standard techniques which utilize the Hardy inequality type spectral
gap estimate to control $|b|v_\pm^2\zeta^2$ in one step.  We choose
$\delta$ sufficiently small in order to absorb this term into the
dissipation.

We have $\partial_r \zeta /r = \partial_\rho \zeta
 /\rho$ since $\zeta$ is radial;  so that 
the singularity $1/\rho$ is effectively $1/R$.  We thus have
\begin{equation}\label{deGiorgiPm}
\begin{aligned}
\sup_{- \sigma^2R^2<t<0} \int_{B_{\sigma R}\times\{t\}}  |(\Omega-k)_\pm|^p  &+ 
\int_{Q_{\sigma R}} |\nabla  (\Omega-k)_\pm^{p/2}  |^2 
\\
&\le 
\frac {C}{(1-\sigma)^2 R^2} \int_{Q_{R}}|(\Omega-k)_\pm|^p.
\end{aligned}
\end{equation}
Our goal will be to establish $L^p$ to $L^\infty$ bounds for functions in this energy class.

\subsection{Local maximum estimate}\label{sec:wH}

The estimates in this section will be proven for a general function $u=\Omega$ satisfying \eqref{deGiorgiPm}:

\begin{lemma}\label{weakH}  Suppose $u=\Omega$ satisfies \eqref{deGiorgiPm} for $1<p\le 2$.  Then
$$
\sup_{Q_{R/2}} u_{\pm} \le C(p, C_*)
\left(R^{-3-2}\int_{Q_{R}}|u_\pm|^p \right)^{1/p}.
$$
\end{lemma}

This estimate can be found in \cite{MR1465184} for $p=2$.  The proof below is similar and we include it so that the proof of Theorem
3.1, which uses Lemma 2.2, is self-contained.  Our choice of $p$ is made merely because those are the ones we need although others are possible.   \\

\begin{proof} 
For $K>0$ to be determined and $N$ a positive integer we define
\begin{gather*}
k_N
=
k_N^\pm=(1\mp 2^{-N})K, 
~R_N=(1+2^{-N})R/2, 
~\rho_N=\frac{R}{2^{N+3}},
\\
R_{N+1}<\bar{R}_N=(R_N+R_{N+1})/2 <R_N.
\end{gather*}
Notice that
$$
R_N-\bar{R}_N
=(R_N-R_{N+1})/2
=(2^{-N}-2^{-N-1})R/4
=\rho_N.
$$
Define $Q_N=Q(R_N)$ and $\bar{Q}_N=Q(\bar{R}_N)\subset Q_N$.
Choose a smooth test function $\zeta_N$ satisfying $\zeta_N\equiv 1$ on $\bar{Q}_N$, $\zeta\equiv 0$ outside $Q_N$ and vanishing on it's spatial boundary, $0\le \zeta_N\le 1$ and 
$\left| \nabla \zeta_N\right|\le \rho_N^{-1}$ in $Q_N$.
Further let
$$
A^\pm(N)=\{X\in Q_N: \pm(u-k_{N+1})(X)>0\}.
$$
And $A_{N,\pm}=\left|A^\pm(N) \right|$.  
Let $v_\pm = \zeta_N (u-k_{N+1})_\pm^{p/2}$.  

H{\"o}lder's inequality gives us
\begin{gather*}
\begin{split}
\int_{Q_{N+1}} ~ |(u-k_{N+1})_\pm|^p
&\le
\int_{\bar{Q}_{N}} ~ |v_\pm|^2
\\
&\le
\left(\int_{\bar{Q}_{N}} ~ |v_\pm|^{2(n+2)/n}\right)^{n/(n+2)}
A_{N,\pm}^{2/(n+2)}.
\end{split}
\end{gather*}
We will use the following parabolic Sobolev inequality:
\begin{equation*}
\int_{Q_R} |u|^{2(n+2)/n}
\le 
C(n) \left( \sup_{-R^2<t<0} \int_{B_R\times\{t\}} |u|^2 \right)^{2/n}
\int_{Q_R}  |\nabla u|^2.
\end{equation*}
See \cite[Theorem 6.11, p.112]{MR1465184}.  We are interested in the form
\begin{equation*}
\int_{Q_R} |u^{p/2}|^{2(n+2)/n}
\le 
C(n) \left( \sup_{-R^2<t<0} \int_{B_R\times\{t\}} |u|^p \right)^{2/n}
\int_{Q_R}  |\nabla u^{p/2}|^2.
\end{equation*}
As in the above followed by Young's inequality then followed by \eqref{deGiorgiPm} we obtain
\begin{gather*}
\left(\int_{\bar{Q}_{N}} ~ |v_\pm|^{2(n+2)/n}\right)^{n/(n+2)}
\\
\le
C \left( \sup_{-R_N^2<t<0} \int_{B(R_N)\times\{t\}} |v_\pm|^2 \right)^{2/(n+2)}
\left(\int_{Q_{N}}  |\nabla v_\pm|^2\right)^{n/(n+2)}
\\
\le
C \left( \sup_{-R_N^2<t<0}\int_{B(R_N)\times\{t\}}|v_\pm|^{2}
+\int_{Q_N}|\nabla v_\pm|^{2}
\right)
\\
\le
C \left( \sup_{-R_N^2<t<0}\int_{B(R_N)\times\{t\}}|(u-k_{N+1})_\pm|^{p}
+\int_{Q_N}|\nabla (u-k_{N+1})_\pm^{p/2}|^{2}
\right)
\\
+\frac{C_*}{\rho_N^2}\int_{Q_N}|(u-k_{N+1})_\pm|^{p}
\\
\le
\frac{C_*}{\rho_N^2}\int_{Q_N}|(u-k_{N+1}) _\pm|^{p}
\le
\frac{C_*}{\rho_N^2}\int_{Q_N}|(u-k_{N}) _\pm|^{p}.
\end{gather*}
Further assume $K^p \ge R^{-n-2}\int_{Q(R)} | u_\pm|^p$.  And define
$$
Y_N \equiv K^{-p} R^{-n-2}\int_{Q_N} | (u-k_{N})_\pm|^p.
$$
Since $k_N^\pm$ are increasing for $+$ or decreasing for $-$ and $Q_N$ are decreasing, $Y_N$ is decreasing.

Chebyshev's inequality tells us that 
\begin{gather*}
A_{N,\pm} =\left|\{ Q_N: \pm(u-k^{\pm}_{N+1})>0  \}\right|
=\left|\{ Q_N: \pm(u-k^{\pm}_{N})>\pm (k^{\pm}_{N+1}- k^{\pm}_N)  \}\right|
\\
=\left|\{ Q_N: \pm(u-k_{N})> K/2^{N+1}  \}\right|
\le  2^{p(N+1)} R^{n+2}Y_N. 
\end{gather*}
Putting all of this together yields
\begin{gather*}
\int_{Q_{N+1}} ~ |(u-k_{N+1})_\pm|^p
\le
\left(\int_{\bar{Q}_{N}} ~ |v_\pm|^{2(n+2)/n}\right)^{n/(n+2)}
A_{N,\pm}^{2/(n+2)},
\\
\le
\left(\frac{C_*}{\rho_N^2}\int_{Q_N}|(u-k_{N}) _\pm|^{p}\right)
\left(2^{p(N+1)} R^{n+2}Y_N\right)^{2/(n+2)}
\\
\le 
\left( \frac{C_*}{\rho_N^2} K^pR^{n+2} Y_N\right)
\left(2^{p(N+1)} R^{n+2}Y_N\right)^{2/(n+2)}
\\
=
C_* K^p
2^{2(N+3)}2^{2p(N+1)/(n+2)} R^{n+2}
Y_N^{1+\frac{2}{n+2}}.
\end{gather*}

We have thus shown that
$$
Y_{N+1}
\le C(N) Y_N^{1+\frac{2}{n+2}}
$$
Here
$
C(N)
=
C_*
2^{2(N+3)}2^{2p(N+1)/(n+2)}.
$
We now choose $K$ as
$$
K^{p}=\left(1+\frac{1}{C_0} \right)R^{-n-2}\int_{Q_0} | u_\pm|^p.
$$
Above
the constant $C_0$ is chosen to ensure that $Y_N\to 0$ as $N\to \infty$.  
\end{proof}

\subsection{Regularity of the original solution}\label{sec:regO}

The limiting solution $\Om$ satisfies \eqref{Om2}, \eqref{Om3} and \eqref{deGiorgiPm}. We conclude from Lemma \ref{weakH} that 
\[
 \Om \in L^\infty(Q_{5/16}).
\]
We further know that $\curl \lv = \lom e_\th \in L^\infty(Q_{5/16})$ from the above estimate on $\Om$ since $\bar{v}_\theta=0$.  Also $\div \lv = 0$ from the equation.
Next $\lv
\in L^\infty_t L^1_x(Q_{5/16})$ by \eqref{assumption}.
We thus conclude $\nb \lv \in L^\infty_t L^4_x(Q_{1/4})$
by Lemma \ref{thA-1}.   Thus $\lv \in L^\infty(Q_{1/4})$ by
embedding.

Now we can deduce regularity of the original solution from
the regularity of the limit solution. 
Since $\lv \in L^\infty(Q_{1/4})$ for $R$ sufficiently small we have
\[
\frac{1}{R^2} \int_{Q_R} |\lv|^3 \le \e_1/2,
\]
where $\e_1$ is the small constant in Lemma \ref{th2-1}.
Fix one such $R>0$. Since $v^{\la_k} \to \lv$ strongly in $L^3$ 
for $k$ sufficiently large we have
\[
\frac{1}{R^2} \int_{Q_R} |v^{\la_k}|^3 \le \e_1.
\]
But this is a contradiction to
\eqref{blowup0}. 
Thus every point $x_*$ on
the $z$-axis is regular; that is, there is a radius $R_{x_*}>0$ so
that $v \in L^\infty(Q(x_*, R_{x_*}))$. Since any finite portion of
the $z$-axis can be covered by a finite subcover of $\{Q(x_*, R_{x_*})\}$,
we have proved Theorem \ref{mainthm}.

The rest of the paper is devoted to proving the key Theorem \ref{keythm}.

\section{H{\"o}lder estimate for axisymmetric solutions}\label{sec:holder}

We now move from cartesian to cylindrical coordinates via the standard
change of variables $x=(x_1, x_2, x_3) =(r\cos\theta, r\sin\theta,
z)$.  For axisymmetric solutions $(v,p)$ of the form \eqref{v-cyl},
the Navier-Stokes equations \eqref{nse} take the form
\begin{align}
\frac{\pd v_r}{\pd t}
+
b\cdot \nb v_r
-
\frac{v_\theta^2 }{r}
+
\frac{\partial p}{\partial r}
&=
\left(\Delta 
-
\frac{1}{r^2}\right) v_r,\nonumber
\\
\frac{\pd v_\theta}{\pd t}
+
b\cdot \nb
v_\theta
+
\frac{v_\theta v_r}{r}
&=
\left(\Delta 
-
\frac{1}{r^2}\right) v_\theta,
\nonumber
\\
\frac{\pd v_z}{\pd t}
+
b\cdot \nb v_z
+
\frac{\partial p}{\partial z}
&=\Delta v_z,\nonumber
\\
 \frac{1}{r} \frac{\partial (r v_r)}{\partial r}
+ \frac{\partial v_z}{\partial z}
 &=0. \nonumber
\end{align}
The vector $b$ is given by
\[
b = v_r e_r + v_z e_z, \quad \div b=0.
\]
The equations of the vorticity $\om=\curl v$, decomposed in the form
\eqref{om-cyl}, are 
\begin{align}
& \frac{\pd \om_r}{\pd t} + b\cdot \nabla \om_r  - \om_r
\pd_r u_r - \om_z \pd_z u_r=\left(\Delta - \frac 1{r^2}\right) \om_r,\nonumber
\\
& \frac{\pd \om_\th}{\pd t} + b\cdot \nabla \om_\th -  2\frac
{u_\th}r \pd_z u_\th - \frac{u_r}r \om_\th = \left(\Delta - \frac 1{r^2}\right) \om_\theta,  \nonumber
\\
& \frac{\pd \om_z}{\pd t}+ b\cdot \nabla \om_z  - \om_z\pd_z u_z  - \om_r
\pd_r u_z=\Delta \om_z .\nonumber
\end{align}
Although we do not use them.
We are interested in the equation for $v_\theta$, which is independent of the pressure.

Consider the change of variable $ \Gamma = r v_\theta$, which is well known (see the references in the introduction). 
The function $\Gamma$ is smooth
and satisfies
\begin{equation}
\frac{\pd \Gamma}{\pd t} + b \cdot \nb \Gamma - \Delta \Gamma +
\frac{2}{r}\frac{\partial \Gamma}{\partial r} = 0.
\label{INSsymC2}
\end{equation}
Note that the sign of the term $\frac{2}{r}\frac{\partial
\Gamma}{\partial r}$ is opposite to that of \eqref{Om2}.  It follows directly from \eqref{assumption} that $\| \Gamma \|_{L^\infty_{t,x}}\le
C_*$; see \cite{MR1902055} for related estimates.
Since $v$ is smooth, we have $\Gamma(t,0,z)=0$ for $t<0$.  The
smoothness and axisymmetry assumptions also imply $v_\th(t,0,z)=0$,
but we will not use this fact.  The main result of this section is the following.

\begin{theorem}\label{keythm}  Suppose that $\Gamma(x,t)$ is a smooth bounded
solution of \eqref{INSsymC2} in $Q_2$ with smooth $b(x,t)$, both may depend on $\th$,
  and
\[
\Gamma|_{r=0}=0, \quad \div b=0, \quad \left| b\right| 
\le C_*/r \quad \text{in }Q_2.
\] 
Then there exist constants $C$ and $\al>0$ which depend only upon
$C_*$ such that
\[
\left| \Gamma(x,t)\right| \le C \| \Gamma \|_{L^\infty_{t,x}(Q_2)} r^\alpha
\quad \text{in }Q_1.
\]
\end{theorem}

We remark that the condition  above is substantially weaker than
\eqref{assumption}, and we do not need $\Gamma$ to be axisymmetric. 
In the rest of this section, we will prove the theorem.
Here we are facing two difficulties: First, the condition $\left. \Gamma
\right|_{r=0}=0$ precludes a direct lower  bound on the fundamental solution and a Harnack inequality on $\Gamma$ (since, when $b=0$, $\Gamma=r^2$ is a nonnegative solution which does not satisfy the usual Harnack inequality.)
Second, the condition $b \le C /r$ is weaker than the standard assumption
$b \le C /|x|$ (see the discussion below). It turns out that one can
develop new techniques incorporating the methods introduced by De
Giorgi \cite{MR0093649} and Moser \cite{MR0170091} to over come
these two points. However, we do not know if one can follow the approach of Nash \cite{MR0100158,MR855753} which relies critically 
on a Gaussian lower bound of the fundamental solution.   The proof of
Theorem 3.1 is independent of the rest of the paper.

The following related  equation has been previously studied by  Zhang \cite{MR2031029}:
$$
\frac{\pd u}{\pd t} + b \cdot \nb u - \Delta u = 0.
$$
He has shown among other things H{\"o}lder continuity of solutions to this equation if $b=b(x)$ is independent of time and $b$ satisfies an integral condition which is fulfilled if say $b$ is controlled by $1/|x|$.  His proof makes use of Moser iteration and Gaussian bounds.


\subsection{Notation, Reformulation, and Energy inequalities}

Let $X=(x,t)$.
Define the modified parabolic cylinder at the origin
$$
Q(R,\tau)=\{X: |x|<R,  - \tau R^2<t< 0\}.
$$
Here $R>0$ and $\tau \in (0,1]$.
We sometimes for brevity write 
$Q_R=Q(R)=Q(R,1)$.
Let
\[
m_2 \equiv \inf_{Q(2R)} \Gamma, 
\quad M_2 \equiv \sup_{Q(2R)} \Gamma,
\quad M \equiv M_2 - m_2>0.
\]
Notice that $m_2 \le 0 \le M_2$ since $\Gamma |_{r=0}=0$.

Now we reformulate the problem in $Q(2R)$ into a new function, $u$,
which will be zero when $\left| \Gamma \right|$ is at its maximum
value.  Specifically, we define
\begin{equation}\label{udef}
u
\equiv
\left\{
\begin{aligned}
2(\Gamma - m_2)/M & \quad \text{if}\quad -m_2 > M_2,
\\
2(M_2 - \Gamma)/M & \quad \text{else.}
\end{aligned}
\right.
\end{equation}
In either case $u$ solves \eqref{INSsymC2} and
$
0 \le u \le 2 
$
in
$Q(2R)$.  We will further  use 
$$
a \equiv u|_{r=0}
=\frac{2}{M}\left(\sup_{Q(2R)}\left| \Gamma \right| \right)
=\frac{2}{M}\max\{M_2, -m_2\}\ge 1,
$$
which follows from our conditions.

We now derive energy estimates for \eqref{INSsymC2}.  Define $v_{\pm}=(u-k)_\pm$ with $k \ge 0$. 
We have $v_+\le (2-k)_+$ and $v_-\le k$.  Consider a  radial test function $0\le \zeta(x,t)\le 1$ for which $\zeta=0$ on 
$\pd B_R \times [-\tau R^2,0]$ and $\frac{\partial \zeta}{\partial r} \le 0$.
We multiply \eqref{INSsymC2} for $u-k$ with $\zeta ^2v_\pm$ and integrate over $\mathbb{R}^3\times [t_0, t]$ to obtain
\begin{equation}\label{eq2}
\begin{aligned}
\frac{1}{2}\bkt{\int_{\mathbb{R}^3}  | \zeta v_\pm|^2  }_{t_0}^t +
&
 \int_{t_0}^t \int_{\mathbb{R}^3} |\nabla(\zeta v_\pm)|^2 
\\ 
=& 
\int_{t_0}^t\int_{\mathbb{R}^3} v_\pm^2 \bke{ b \zeta \cdot
\nabla \zeta + \zeta \frac{\partial\zeta}{\partial t} + |\nabla \zeta|^2 + \frac {2 \zeta}r
\frac{\partial\zeta}{\partial r} } 
\\
&+ 2\pi [(a-k)_\pm]^2\int_{t_0}^t\int_{\mathbb{R}}dz ~\zeta^2|_{r=0}.
\end{aligned}
\end{equation}
We need to estimate all the terms in parenthesis.

Choose $\sigma\in(1/4,1)$, we require that the test function satisfies $\zeta \equiv 1$ on $Q(\sigma R,\tau)$.
If we further choose $\zeta(x,t_0)=0$ then, using \eqref{b-term}, we  estimate \eqref{eq2} as follows
\begin{equation}\label{eq3}
\begin{aligned}
\sup_{-\tau \si^2R^2<t<0} \int_{B(\si R)\times\{t\}}  v^2_\pm  &+ 
\int_{Q(\si R,\tau)} |\nabla v_\pm|^2 
\\
&\le 
\frac {C_{**}}{\tau
(1-\sigma)^2R^2} \int_{Q(R,\tau)}v_\pm^2 + C \tau R^3 [(a-k)_\pm]^2.
\end{aligned}
\end{equation}
If we alternatively choose $\zeta = \zeta(x)$ then \eqref{eq2} takes the form
\begin{equation}\label{eq4}
\begin{aligned}
\sup_{t_0<s<t} \int_{B(\si R)\times\{s\}}  v^2_\pm  &+ 
\int_{t_0}^t\int_{B(\si R)} |\nabla v_\pm|^2 
-\int_{B_R\times\{t_0\}}v_\pm^2
\\
&\le 
\frac {C_{**}}{(1-\sigma)^2R^2} 
\int_{t_0}^t\int_{B_R}v_\pm^2 + C \tau R^3 [(a-k)_\pm]^2.
\end{aligned}
\end{equation}
Notice that there is no $\tau^{-1}$ appearing in this energy
inequality \eqref{eq4} compared to \eqref{eq3}.

The energy estimates \eqref{eq3} and \eqref{eq4} are the standard
parabolic De Giorgi classes except for the last term.  Our goal will be to use them to show that
the set where $\Gamma$ is very close to its largest absolute value
or, equivalently, the set where $u$ is almost zero is as small as you
wish.  We establish this fact in the following series of Lemma's.
 
 \subsection{Initial Estimates}
 
Later on we will use the two standard Lemma's below in a non-standard
iteration scheme of sorts to show that the set where $u$ is almost
zero has very small Lebesgue measure.
 
 \begin{lemma}\label{timecontinuity}  
Suppose there exists a $t_0\in [-\tau R^2, 0]$, $K>0$ and $\gamma\in (0,1)$ so that
 \begin{equation*} 
\left| \left\{x \in B_R: u(x,t_0) \le K\right\}\right|
\le \gamma |B_R|.
\end{equation*}
Further suppose that  $u$ satisfies \eqref{eq4} for $v_-$.
Then for all $\eta\in (0, 1-\sqrt{\gamma})$ and $\mu\in (\gamma/(1-\eta)^2, 1)$
there exists $\theta \in (0,1)$ such that
 \begin{equation*} 
 \left| \left\{x \in B_R: u(x,t) \le \eta K\right\}\right|
\le \mu |B_R|, \quad \forall t\in [t_0, t_0+(\tau\wedge \theta) R^2].
\end{equation*}
Here $\theta$ depends only on the constants in \eqref{eq4} and $\gamma$.
\end{lemma}

We note that the proof shows that $\theta(\gamma)\to 0$ as $\gamma\uparrow 1$, but if $\tau$ is sufficiently small then we may take 
$\theta=\tau$ when $\gamma$ is close enough to zero.  And if $\gamma$ is small, then $\mu$ can be taken almost as small.  \\

\begin{proof}  We consider $v_-=(u-K)_-$ and assume without loss of generality that $K<1$. 
The energy inequality \eqref{eq4}  for this function is
\begin{align*}
 \int_{B(\si R)\times\{t\}}  v^2_-  
 &\le 
 \int_{B_R\times\{t_0\}}v_-^2 +
\frac {C_{**}}{(1-\sigma)^2R^2} \int_{t_0}^t\int_{B_R}v_-^2 + C
\tau R^3 [(a-K)_-]^2
\\
&\le  K^2 |B_R| \left(\gamma + \frac {C_{**}(\tau\wedge \theta)}{(1-\sigma)^2}\right).
\end{align*}
We have used $(a-K)_-=0$.
The Chebyshev inequality tells us that
\begin{equation*}
A\equiv |\bket{x \in B(\si R): u(x,t) \le \eta K}|\cdot (K-\eta K)^2
\le 
 \int_{B(\si R)\times\{t\}}  v^2_-.
\end{equation*}
The quantity $1-\sigma^3$ is an upper bound for the measure of $A^c$, which grants the following general inequality
\begin{equation*}
\begin{aligned}
\frac{|\bket{x \in B_R, u(x,t) \le \eta K}| }{|B_R|} 
\le 
\frac{|\bket{x \in B(\sigma R), u(x,t) \le \eta K}|}{|B_R|}  + (1-\sigma^3),
\\
\le 
(1-\eta)^{-2}\left(\gamma + \frac {C_{**}(\tau\wedge \theta)}{(1-\sigma)^2}\right)
+(1-\si^3).
\end{aligned}
\end{equation*}
Now let $\sigma$ be so close to one that  
$
\frac{\gamma}{(1-\eta)^{2}}+(1-\si^3) < \mu.
$
Then, with $\tau$ fixed, choose $\theta$ small enough that the whole thing is $\le \mu$.  
\end{proof}

The Lemma above shows continuity in time of the Lebesgue measure of the set where $u$ is small and the lemma below shows that if the set where $u$ is small is less than the whole set, then the set where $u$ is even smaller can be made tiny.  This is an extremely weak way to measure diffusion.

\begin{lemma}\label{smallset}  Suppose that $u(x, t)$ satisfies \eqref{eq3} for $v_-$.  In addition 
\begin{equation*} 
\left| \left\{x \in B_R: u(x,t) \le K\right\}\right|
\le \gamma |B_R|, 
\quad
\forall t\in [t_0, t_0+\theta R^2]=I,
\end{equation*}
where $K,\theta>0$, $\gamma\in (0,1)$ and $B_R\times I \subset Q(R,\tau)$.  Then for all
$\epsilon \in (0,1)$ there exists a $\delta\in(0,1)$ such that
\begin{equation*} 
\left| \left\{X \in B_R\times I: u(X) \le \delta  \right\}\right|
\le \epsilon |B_R\times I|.
\end{equation*}
\end{lemma}

\begin{proof} 
We denote, for $n=0,1,2,3,\ldots$,
$$
A_n(t) =\bket{x \in B_R: u(x,t) \le 2^{-n}K}, 
\quad 
A_n = \cup _{t\in I} A_n(t).
$$
Clearly $|A_{n+1}|\le |A_{n}|\le |A_{0}|\le \gamma |B_R\times I|$.  
And
$$
|A_{n}^c(t)|=\left|\left\{x \in B_R: u(x,t) > 2^{-n}K\right\}\right|
=|B_R|-|A_{n}(t)|\ge (1-\gamma) |B_R|.
$$
Since $\gamma<1$, we know that $A_{n}^c(t)$ does not have measure zero.

We invoke the following well known version of the Poincar{\'e} inequality.  For any $v\in W^{1,1}(B_R)$ and for any $\alpha, \beta\in \mathbb{R}$ with $\alpha<\beta$ we have
$$
|\{x\in B_R: v(x)\le \alpha\}|  \le 
\frac{C R^{3+1}/(\beta-\alpha)}{|\{x\in B_R: v(x)> \beta\}|} \int_{B_R\cap \{\alpha<u\le \beta\}} |\nabla v|,
$$
where $C>0$ only depends on the dimension.  Let $\beta=2^{-n}K$ and $\alpha=2^{-n-1}K$.  We have
$$
|A_{n+1}(t)| 
\le 
\frac{C 2^{n+1}R}{K(1-\gamma)} 
\int_{A_n(t)-A_{n+1}(t)} |\nabla u|
=
\frac{C 2^{n+1}R}{K(1-\gamma)} 
\int_{A_n(t)-A_{n+1}(t)} |\nabla (u-a)^-|.
$$
We use the Cauchy-Schwartz inequality to bound this integral as
\begin{align*}
|A_{n+1}| 
=
\int_{I} |A_{n+1}(t)| 
&\le 
\frac{C 2^{n+1}R}{K(1-\gamma)} 
\int_{A_n-A_{n+1}} |\nabla (u-a)^-|
\\
&\le 
\frac{C 2^{n+1}R}{K(1-\gamma)} 
 |A_n - A_{n+1}|^{1/2} 
 \bke{ \int_{A_n-A_{n+1}} |\nabla (u-a)^-|^2}^{1/2}.
\end{align*}
The energy inequality \eqref{eq3}, with $\si R$ and $R$ replaced by $R$ and $2R$ results in
\begin{align*}
|A_{n+1}| 
&\le
\frac{C 2^{n+1}R}{K(1-\gamma)} 
 |A_n - A_{n+1}|^{1/2} 
 \bke{\frac {C}{\tau R^2} \int_{Q(2R, \tau)} |(u-a)^-|^2}^{1/2}
\\
& \le 
\frac{C 2^{n+1}R}{K(1-\gamma)}  |A_n - A_{n+1}|^{1/2} |B(2R)|^{1/2}a
=
\frac{C R^{5/2}}{K(1-\gamma)}  |A_n - A_{n+1}|^{1/2}.
\end{align*}
Square both sides of this inequality and dividing by  
$|B_R \times I|^2$ 
 to obtain
$$
\frac{|A_{n+1}|^2}{|B_R \times I|^2} 
\le  
\frac{C}{\theta K^2(1-\gamma)^2}  
\left(\frac{|A_{n}|}{|B_R \times I|}  - \frac{|A_{n+1}|}{|B_R \times I|} \right).
$$
Summing in $n$, we get 
\begin{equation*}
\begin{aligned}
n \frac{|A_{n}|^2}{|B_R \times I|^2}  
\le 
\sum_{j=1}^n \frac{|A_{j}|^2}{|B_R \times I|^2} 
&\le 
\frac{C}{\theta K^2(1-\gamma)^2}   
\sum_{j=1}^n 
\left(\frac{|A_{j-1}|}{|B_R \times I|}  - \frac{|A_{j}|}{|B_R \times I|} \right)
\\
& =
\frac{C}{\theta K^2(1-\gamma)^2}   
\left(\frac{|A_{0}|}{|B_R \times I|}  - \frac{|A_{n}|}{|B_R \times I|} \right)
\\
&  \le 
\frac{C}{\theta K^2(1-\gamma)^2}   
\frac{|A_{0}|}{|B_R \times I|}
\le
\frac{C\gamma}{\theta K^2(1-\gamma)^2}   .
\end{aligned}
\end{equation*}
We complete the proof by choosing $n$ sufficiently large.
\end{proof}
 

 \subsection{Estimate on the measure of the set where $u$ is small}
The next lemma allows us to apply all the machinery above.

\begin{lemma}\label{eq5}
There exists a $\kappa\in (0,1)$ such that
$0<\la < \min\{\kappa \tau, 1/8\}$ implies
\begin{equation*}
 \left|\left\{X \in Q(R,\tau): u(X) \le \la^2\right\}\right| \le (1-4\la) |Q(R,\tau)|.
\end{equation*}
\end{lemma}

\begin{proof}We establish a contradiction using energy estimates.
Suppose the opposite
$$
 \left|\left\{X \in Q(R,\tau): u(X) \le \la^2\right\}\right| >(1- 4\la) |Q(R,\tau)|.
$$
Or equivalently
\begin{equation}\label{eq5-}
 \left|\left\{X \in Q(R,\tau): u(X) > \la^2\right\}\right| < 4\la |Q(R,\tau)|.
\end{equation}
This condition will imply a contradiction to the size condition on $a\ge 1$.  

We will test the equation \eqref{INSsymC2} with 
$pu^{p-1} \zeta^2$ for $0<p<1$ and $\zeta \ge 0$.  
Since $u=0$ sometimes, in general we should test \eqref{INSsymC2} for $u+\epsilon$ with $p(u+\epsilon)^{p-1} \zeta^2$ and then send $\epsilon \downarrow 0$ to obtain our estimates.  However, since the result is the same, to simplify the presentation we will omit these details.  We have
\begin{equation*}
\begin{aligned}
\int_{Q(R,\tau)} pu^{p-1} \zeta^2 \frac{\pd u}{\partial t}  
= &
\bkt{\int_{B_R} \zeta^2 u^p}^0_{t_1} - \int_{Q(R,\tau)} u^p 2 \zeta \frac{\pd \zeta}{\partial t}
\equiv I_1 + I_2,
\\
\int_{Q(R,\tau)} pu^{p-1} \zeta^2 (-\Delta u) 
=& \frac {4(p-1)}p \int_{Q(R,\tau)} |\nabla
(u^{p/2}\zeta)|^2 
\\
&+ \int_{Q(R,\tau)} 2 u^p \bkt{-|\nabla \zeta|^2 + \frac
{p-2}p \zeta\Delta \zeta}
\equiv I_3 + I_4,
\\
\int_{Q(R,\tau)}  pu^{p-1} \zeta^2 b \cdot \nabla u 
= &
- \int_{Q(R,\tau)}  2u^p b  \cdot \zeta \nabla \zeta 
\equiv I_5,
\\
\int_{Q(R,\tau)} pu^{p-1} \zeta^2 \frac 2r \pd_r u 
= & - \int_{Q(R,\tau)} 4 u^p\zeta \zeta _\rho
/\rho - \int_{-\tau R^2}^{0} dt\int_{\mathbb{R}} dz ~ 2 (\zeta^2 u^p)|_{r=0} 
\\
\equiv I_6 + I_7.
\end{aligned}
\end{equation*}
In the computation of $I_6$ we have used $\zeta _r /r = \zeta
_\rho /\rho$, which follows if $\zeta=\zeta(\rho, t)$ where 
$\rho=|x|=\sqrt{r^2+z^2}$. Notice that $\sum_{j=1}^7 I_j=0$.
For arbitrary $p\in (0,1)$, we see that $I_3$ and $I_7$ are both
non-positive.

We choose $\zeta = \zeta_1(\rho) \zeta_2(t)$ where
$\zeta_1(\rho) =1$ in $B(R/2)$ and $\zeta_1(\rho)$ has compact support
in $B_R$; also $\zeta_2(t) =1$ if $t \in [-\frac 78 \tau R^2,
-\frac 18 \tau R^2]$ and $\zeta_2(t)$ has compact support in
$(-\tau R^2,0)$. Thus $I_1=0$ and we have
\[
\frac{6}{4}\tau R^3 a^p \le -I_7 =\sum_{j=2}^6 I_j.
\]
We estimate each of the terms $I_2$ through $I_6$ to obtain a contradiction.

By the argument in \eqref{b-term}, we have
$$
|I_5| \le \frac{2(1-p)}{p} \int_{Q(R,\tau)} |\nabla
(u^{p/2}\zeta)|^2
 + \frac{C}{R^2} \int_{Q(R,\tau)} u^p. 
$$
Also note $\nabla \zeta=0$ in $B(R/2)$ and so
the singularity $1/\rho$ is effectively $1/R$. Thus,
\[
I_2 \le \frac {C}{\tau R^2} \int_{Q(R,\tau)} u^p, \quad \sum_{j=3}^6
I_j \le \frac {C}{R^2} \int_{Q(R,\tau)} u^p .
\]
Assuming \eqref{eq5-} and using $0\le u \le 2$, we have
\[
a^p \le \frac {C}{\tau^2 R^5} \int_{Q(R,\tau)} u^p \le \frac
{C}{\tau^2 R^5}
\left\{ \la^{2p} |Q(R,\tau)| + 2^p ( 4 \la |Q(R,\tau)|) \right\}
\le \frac {C_2}\tau (\la^{2p} + \la).
\]
Here $C_2=C_2(C_*)$. 
Take $p=1/2$ and $\kappa = \frac{1}{4C_2}$ to get $a^p < 1$, a contradiction.\end{proof}

\bigskip

Lemma \ref{eq5} is the starting point of our iteration scheme.
From this Lemma we know that there is a 
$t_1 \in [-\tau R^2, - 2\la \tau R^2]$ 
so that
\begin{equation}\label{eq6}
|\bket{x \in B_R: u(x,t_1) \le \la^2}| 
\le (1-2\la) |B_R|.
\end{equation}
Then apply Lemma \ref{timecontinuity} with $K=\la^2$ to \eqref{eq6} to see, 
for say $\eta=\la$ and $\mu=1-\la$, that
\begin{equation*}
|\bket{x \in B_R: u(x,t) \le \la^3}| \le (1-\la) |B_R|,
\quad \forall t\in[t_1, t_1+\theta_* R^2]\equiv I_*.
\end{equation*}
Here $\theta_*=\theta \wedge \tau$ and $\theta$ is the constant chosen in Lemma \ref{timecontinuity}. 
From here Lemma \ref{smallset} allows us to conclude
\begin{equation*}
|\bket{X \in B_R\times I_*: u(X) \le \delta_*}| \le \frac{\epsilon_*}{2} |B_R\times I_*|,
\end{equation*}
where $\epsilon_* >0$ is as small as you want and 
$\delta_*=\delta_*(\epsilon_*)$.

Then, as in \eqref{eq6}, there exists a $ t_2\in I_*$ (so that $t_2\le -\la \tau R^2$) such that
\begin{equation}
|\bket{x \in B_R: u(x, t_2) \le \delta_*}| 
\le 
\epsilon_* |B_R|.
\label{epsilon}
\end{equation}
Up till now all the small parameters that we have chosen depend upon $\tau$. 
But above $\epsilon_*$  can be taken arbitrarily small independent of the size of $\tau$.  This is the key point that enables us to proceed.  It is the reason why we are required to do this procedure twice.

Now suppose 
$1-\sigma^3=1/4$ and choose first $\tau<1/8$  so that 
$C_{**}\tau/(1-\sigma)^2\le 1/4$.  
Then take $\delta_*$ from \eqref{epsilon} with $\epsilon_* <1/16$ above playing the role of $\gamma$ in Lemma \ref{timecontinuity}.  Also $\eta< 1/2$.  With all this, from Lemma \ref{timecontinuity}, we can choose $\mu<1$ so that 
\begin{equation*}
\left|\left\{x \in B_R: u(x,t) \le \eta \delta_*\right\}\right| \le \mu |B_R|,
\quad \forall t\in[t_2, t_2+\tau R^2]\equiv I.
\end{equation*}
Further, it is safe to assume 
 that $\theta_* \le \la $;  we see that 
 $t_2\le -\la\tau R^2$ and so $[-\la\tau R^2, 0]\subset I$.
Finally apply Lemma \ref{smallset} again to obtain
\begin{equation}\label{keysmallness}
\left|\left\{X \in Q(R,\la\tau): u(X) \le \delta\right\}\right|
\le 
\epsilon |Q(R,\la\tau)|,
\end{equation}
with $\epsilon>0$ arbitrarily small.  This is a key step in what follows.

Let $U=\delta -u$, where $\delta$ is the constant from \eqref{keysmallness}. $U$ is clearly a solution of
\eqref{INSsymC2} and $U|_{r=0} =\delta -a< 0$.  We apply \eqref{eq3} to $U$ on
$Q(2d)$ (with $\tau =1$) to get
\begin{equation*}
\begin{aligned}
\sup_{-\si^2d^2<t<0} \int_{B(\si d)\times\{t\}}  |(U-k)^+|^2  &+ 
\int_{Q(\si d)} |\nabla (U-k)^+|^2 
\\
&\le 
\frac {C_{**}}{
(1-\sigma)^2d^2} \int_{Q(d)}|(U-k)^+|^2.
\end{aligned}
\end{equation*}
This holds for all $k>0$ and $\si \in (0,1)$. So we can
apply Lemma \ref{weakH} to conclude
\begin{equation}\label{eq15}
 \sup_{Q(d/2)} ( \delta - u)  \le \bke{ \frac{C}{|Q(d)|}
\int _{Q(d)} |(\delta -u)^+|^2}^{1/2}.
\end{equation}
This  inequality combined with \eqref{keysmallness} will produce a lower bound.

\subsection{Regularity from a lower bound}

Let $d= \sqrt{\la \tau } R$
so that
$Q(d) \subset Q(R,\la \tau)$.  By \eqref{eq15}
and \eqref{keysmallness},
\begin{align*}
\delta - \inf_{Q(d/2)} u & \le \bke{ \frac{C}{|Q(d)|} \int
_{Q(d)} |(\delta -u)^+|^2}^{1/2}
\\
& \le 
\left( \frac{C\delta^2 \epsilon |Q(R, \la \tau)|}{|Q(d)|}  
\right)^{1/2}
=
C \delta \epsilon^{1/2} \left( \la \tau\right)^{-3/4},
\end{align*}
which is less than $\frac \delta 2$ if $\epsilon$ is chosen sufficiently
small. We conclude
\[
\inf_{Q(d/2)} u \ge \frac \delta 2 .
\]
This is the lower bound we seek.  From it we will deduce an oscillation estimate.

This entails a bit of algebra. We define
\[
m_d \equiv \inf_{Q(d/2)} \Gamma, 
\quad M_d \equiv \sup_{Q(d/2)} \Gamma.
\]
Then from \eqref{udef} we have
$$
\inf_{Q(d/2)} u =
\left\{
\begin{aligned}
2(m_d - m_2)/M & \quad \text{if}\quad -m_2 > M_2,
\\
2(M_2 - M_d)/M & \quad \text{else},
\end{aligned}
\right.
$$
Notice that both expressions above are non-negative in any case; thus we can  
add them together to observe that
$$
\frac \delta 2 
\le
\frac{2}{M}\left\{M-\osc(\Gamma,d/2)\right\}
$$
Here 
$\osc(\Gamma,d/2)=M_d-m_d$ 
and 
$\osc(\Gamma,2R)=M_2-m_2=M$.  
We rearrange the above 
$$
\osc(\Gamma, d/2) \le 
\left(1-\frac{\delta}{4} \right) 
\osc(\Gamma, 2R).
$$
This is enough to produce the desired H{\"o}lder continuity via the following.

\subsection{Iteration Argument}
Suppose we have a non-decreasing function $\omega$ on an interval $(0, R_0]$ which satisfies
$$
\omega(\tau R) \le \gamma \omega(R),
$$
with $0<\gamma, \tau<1$.  Then for $R\le R_0$ we have
\begin{equation}\label{iteration}
\omega(R)\le \frac{1}{\gamma} \left( \frac{R}{R_0}\right)^\alpha \omega(R_0),
\end{equation}
where  $\alpha=\log\gamma/\log\tau>0$.

Iterating, as in \eqref{iteration}, we get, for  
$C_\Gamma=\left(1-\frac{\delta}{4} \right)^{-1} \sup_{Q(1)} \Gamma$, that
\begin{equation}\label{eq20}
 \osc(\Gamma, R) \le  C_\Gamma 
 R^\al,
 \quad \forall R \in (0,1),
\end{equation}
for 
$\al =2 \log(1-\frac{\delta}{4})/\log(\la \tau/16)>0$. 
Thus $\Gamma$ is H{\"o}lder continuous near the origin.
We have proved Theorem \ref{keythm}.

\section*{Appendix}
\addcontentsline{toc}{section}{Appendix}
\setcounter{theorem}{0}
\renewcommand{\thetheorem}{A.\arabic{theorem}}
\setcounter{equation}{0}
\renewcommand{\theequation}{A.\arabic{equation}}

Here we collect some estimates needed for Section 2.

\begin{lemma}\label{thA-1}
Let $B_{R_2} \subset B_{R_1} \subset \R^3$ be concentric with $0<R_2<R_1$.  Let $v$
be a vector field defined in $B_{R_1}$. Let $1<q<\infty$ and $0<\al<1$. Then for
$k=0,1,\ldots$ there is a constant $c$ depending on ${R_2},{R_1},q,\al,k$ so that
\begin{equation*}
\norm{\nb^{k+1} u}_{L^q(B_{R_2})} \le c \norm{\nb^k \div u}_{L^q(B_{R_1})}
+ c \norm{\nb^k \curl u}_{L^q(B_{R_1})} + c \norm{ u}_{L^1(B_{R_1})}.
\end{equation*}
and
\begin{equation*}
\norm{\nb^{k+1} u}_{C^\al(B_{R_2})} \le c \norm{\nb^k \div u}_{C^\al(B_{R_1})}
+ c \norm{\nb^k \curl u}_{C^\al(B_{R_1})} + c \norm{ u}_{L^1(B_{R_1})}.
\end{equation*}
\end{lemma}

This is well-known, see \cite{MR1397564}.

\begin{lemma}[Interior estimates for Stokes system]\label{thA-2}
Fix $R \in (0,1)$. Let $1<s,q<\infty$ and $f=(f_{ij}) \in
L^s_tL^q_x(Q_1)$. Assume that $v \in L^s_tL^1_x(Q_1)$ is a weak solution of
the Stokes system
\[
\pd_t v_i - \Delta v_i + \pd_i p = \pd_j f_{ij}, \quad \div v = 0
\quad \text{in }Q_1.
\]
Then $v$ satisfies, for some constant $c=c(q,s,R)$,
\begin{equation} \label{eqA1}
\norm{\nb v}_{L^s_tL^q_x(Q_R)} \le c \norm{f}_{L^s_tL^q_x(Q_1)} + c
\norm{v}_{L_t^sL_x^1(Q_1)}.
\end{equation}
If instead $v$ is a weak solution of
\[
\pd_t v_i - \Delta v_i + \pd_i p = g_{i}, \quad \div v = 0
\quad \text{in }Q_1,
\]
then
\begin{equation}\label{eqA2}
\norm{\nb^2 v}_{L^s_tL^q_x(Q_R)} \le c \norm{g}_{L^s_tL^q_x(Q_1)} + c
\norm{v}_{L^sL^1(Q_1)}.
\end{equation}
\end{lemma}

An important feature of these estimates is that a bound of the
pressure $p$ is not needed in the right side. A similar estimate for
the time-independent Stokes system appeared in \cite{MR1789922}.
Note that these estimates improve the spatial regularity only. One
cannot improve the temporal regularity,
in view of Serrin's example of a solution $v(x,t) = f(t) \nb h(x)$ where
$h(x)$ is harmonic.

\medskip

\begin{proof} \quad 
Denote by $P$ the Helmholtz projection in $\R^3$, $(Pg)_i = g_i - R_i
R_k g_k$, where $R_i$ is the $i$-th Riesz transform.  Let $\tau =
R^{1/4} \in (R,1)$ and choose $\zeta(x,t)\in C^\infty (\R^4)$,
$\zeta \ge 0$, $\zeta = 1$ on $Q_\tau$ and $\zeta=0$ on $\R^3 \times
(-\infty,0] - Q_1$.  For a fixed $i$, define
\[
\td v_i(x,t) =
\int_{-1}^t \Gamma(x-y,t-s) \, \pd_j ( F_{ij})(y,s)\,dy\,ds,
\]
where $\Gamma$ is the heat kernel and $ F_{ij} = f_{ij} \zeta - R_i R_k
(f_{kj}\zeta)$. The function $\td v_i$ satisfies 
$$
(\pd_t - \Delta)\td v_i = \pd_j F_{ij} = [P\pd_j\zeta( f_{kj})_{k=1}^3]_i, 
\quad
\div \td v = 0.
$$ 
The $L^s_tL^q_x$-estimates for the parabolic version of singular
integrals and potentials (see \cite{MR0234130,MR0440268,MR1232192}, also see \cite{MR657505}, \cite{MR1865493} and their references), and the usual version of $L^q$-estimates
for singular integrals,   give
\begin{equation}\label{eqA3}
\norm{\nb \td v}_{L^s_tL^q_x(Q_1)} +\norm{ \td v}_{L^s_tL^q_x(Q_1)}\le 
c\norm{F}_{L^s_t L^q_x} \le c \norm{f}_{L^s_t L^q_x(Q_1)}.
\end{equation}
Furthermore, for some function $\td p(x,t)$,
\[
(\pd_t - \Delta) \td v + \nb \td p = \pd_j(\zeta f_{ij}), \quad \div \td
v = 0.
\]
The differences $u =v - \td v$ and $\pi = p -\td p$ satisfy the
homogeneous Stokes system
\[
\pd_t u - \Delta u + \nb \pi = 0, \quad \div v = 0
\quad \text{in }Q_\tau.
\]
Its vorticity $\om = \curl u$ satisfies the heat equation $(\pd_t -
\Delta)\om=0$. Let $W = \zeta_\tau \om$, where $\zeta_\tau(x,t)=
\zeta(x/\tau,t/\tau^2)$.  It satisfies
\[
(\pd_t - \Delta)W = G :=w(\pd_t-\Delta) \zeta_\tau - 2 (\pd_m \zeta_\tau)
\pd_m \om.
\]
And thus, for $(x,t)\in Q_{\tau^2}$,
\[
\om_i(x,t)= W_i(x,t) = \int_{-1}^t \int \Gamma(x-y,t-s) \, G_i(y,s) \,dy ds =
\int_{-1}^t \int H_{x,y}^{i,j} (y,s) \, u_j(y,s) \,dy ds
\]
where, using $\om_i = - \de_{ijk} \pd_k u_j$,
\[
H_{x,y}^{i,j} (y,s) = \pd_{y_k} \de_{ijk} \bket{
\Gamma(x-y,t-s) (\pd_t-\Delta) \zeta_\tau + 2 \div
[\Gamma(x-y,t-s) \nb \zeta_\tau]}.
\]
The functions $H_{x,y}^{i,j}$ are smooth with uniform $L^\infty$-bound
for $(x,t)\in Q_{\tau^3}$. Thus
\[
\norm{\curl u}_{L^\infty(Q_{\tau^3})} \le C \norm{u}_{L^1(Q_1)}.
\]
Since $\div u=0$, we have for any $q<\infty$, using Lemma \ref{thA-1},
\begin{equation}\label{eqA4}
\norm{\nb u}_{L^s_t L^q_x(Q_R)} \le c\norm{u}_{L^s_tL^1_x(Q_1)} 
\le c \norm{v}_{L^s_t L^1_x(Q)} + c\norm{\td v}_{L^s_t L^1_x(Q)}.
\end{equation}
The sum of \eqref{eqA3} and \eqref{eqA4} gives \eqref{eqA1}. The proof
of \eqref{eqA2} is similar: one defines
\[
\td v_i(x,t) =
\int_{-1}^t \Gamma(x-y,t-s) \,   F_i(y,s)\,dy\,ds, \quad
F_i = g_i \zeta - R_i R_k(g_k \zeta)
\]
and obtains $\norm{\nb^2 \td v}_{L^s_t L^q_x(Q_1)} +\norm{ \td
v}_{L^s_t L^q_x (Q_1)}\le c\norm{F}_{L^s_t L^q_x} \le c
\norm{g}_{L^sL^q(Q_1)}$. One then estimates $\norm{\nb^2 (v-\td
v)}_{L^s_t L^q_x (Q_R)}$ in the same way.
\end{proof}


\section*{Acknowledgments}

The authors would like to thank the National Center for Theoretical Sciences (NCTS), Taipei Office, and National Taiwan University for hosting part of our collaboration in the summer of 2006. Tsai would also like to thank Harvard University for their hospitality during the Fall of 2006. The research of Chen is partly supported by the NSC grant 95-2115-M-002-008 (Taiwan). The research of Strain is partly supported by the NSF fellowship DMS-0602513 (USA). The research of Yau is partly supported by the NSF grant DMS-0602038. (USA). The
research of Tsai is partly supported by an NSERC grant (Canada).


\end{document}